\documentstyle[12pt]{article}

\begin{document}
\thispagestyle{empty}
\baselineskip=20pt      
\setcounter{page}{0}
\ \ \ 
\vspace{2.0cm}
\begin{center}
{\large \bf  Fermion Representation of  Quantum Group  \\}
\vspace{1.5cm}
Minoru HIRAYAMA and Shiori KAMIBAYASHI \\
\vspace{1.5cm}
{\it Department of Physics, Toyama University, Toyama 930}\\
\vspace{0.5cm}
{\small (Received November 28, 1996)}
\end{center}
\vspace{3cm}
                             
{\small The spinor representation of the quantum group $U_q(su(N))$ is given in terms of a set of fermion creation and annihilation operators. It is shown that the $q$-fermion operators introduced earlier can be identified with the conventional fermion operators. Algebra homomorphisms mapping the fermion operators to their tensor products are discussed. The relation of the coproduct of the quantum group to the above algebra homomorphisms is obtained. }

\newpage
\begin{center}
{\large \bf \S 1.  Introduction}\par
\end{center}
\vspace{0.5cm}
\ \ \ Quantum groups are  subject much discussed both by physicists and mathematicians. A quantum group  is mathematically defined as a noncommutative and nonco-commutative Hopf algebra.$^{1)-3)}$ The representation theory of quantum groups has also been fully developed and has found many useful applications in physics.$^{4)}$  It is interesting that the analogue of the Jordan-Schwinger representation of Lie groups exists for quantum groups. Biedenharn$^{5)}$ and Macfarlane$^{6)}$ introduced $q$-bosonic operators satisfying $q$-deformed commutation relations and constructed the $q$-boson representation of the simplest case of quantum groups. A thorough discussion on the $q$-boson and $q$-fermion representations of quantum groups was given by Hayashi.$^{7)}$ His discussion covers   various quantum groups $U_q(X)$ with $X=A_{N-1}, B_N, C_N, D_N$ and $A_{N-1}^{(1)}$. We, however,  concentrate here on the quantum group $U_{q}(A_{N-1})$ and assume that the parameter $q$ is real and!
 not  equal to $\pm1$.
The quantum group $U_q(A_{N-1})$ is the enveloping algebra generated by $e_i, f_i,k_i, k_i^{-1} (1\leq i \leq N-1, 2\leq N)$ satisfying the follwing relations:
$$
k_ik_i^{-1}=k_i^{-1}k_i=1, \ \ \ k_ik_j=k_jk_i,   \eqno(1\cdot 1{\rm a}) $$
$$k_ie_jk_i^{-1}=q^{a_{ij}}e_j, \ \ \ k_if_jk_i^{-1}=q^{-a_{ij}}f_j,  \eqno(1\cdot 1{\rm b}) $$
$$
e_if_j-f_je_i=\delta_{ij}\frac{k_i^2-k_i^{-2}}{q^2-q^{-2}},   \eqno(1\cdot 1{\rm c})$$

$$
{\displaystyle\sum_{0\leq n\leq{1-a_{ij}}}}(-1)^n {1-a_{ij}\brack n}_{q^2}e_i^{1-a_{ij}-n}e_je_i^n=0,    \ \ \      (i\neq j)  \eqno(1\cdot1{\rm d}) $$
$$
\displaystyle{\sum_{0\leq n\leq{1-a_{ij}}}}(-1)^n  {1-a_{ij}\brack n}_{q^2}f_i^{1-a_{ij}-n}f_jf_i^n=0. \ \ \         (i\neq j) \eqno(1\cdot 1{\rm e}) $$
In $(1\cdot 1)$, $a_{ij}$ is the $ij-$element of the $(N-1)\times(N-1)$ generalized Cartan matrix
$$
A_{N-1}=\left[\begin{array}{ccccccc}
2 & -1 & 0 & . & . & . & 0 \\  
-1 & 2 & -1 & . & . & . & . \\
0 & -1 & 2 & . & . & . & . \\
. & . & . & . & . & . & .\\
. & . & . & . & 2 & -1 & 0 \\
. & . & . & . & -1 & 2 & -1 \\
0 & . & . & . & 0 &-1 & 2 \end{array} \right].  \eqno(1\cdot2) $$
In $(1\cdot 1)$, we have made use of the following notation: 
$$  {m\brack n}_Q=\frac {[m]!_Q}{[m-n]!_Q[n]!_Q}\ ,$$
$$ [m]!_Q=[m]_Q[m-1]_Q{\cdots}[2]_Q[1]_Q \ ,    $$
$$ [m]_Q=\frac {Q^m-Q^{-m}}{Q-Q^{-1}} \ .\eqno(1\cdot3)  $$
Hayashi$^{7)}$ showed that the relations $(1\cdot1)$ are realized by putting
$$ e_i=\psi_i\psi_{i+1}^\dagger \ ,\eqno(1cdot4\rm a) $$
$$ f_i= \psi_{i+1}\psi_i^\dagger \ , \eqno(1\cdot4\rm b) $$
$$ k_i=\omega_i\omega_{i+1}^{-1},\ \ \  k_i^{-1}=\omega_{i+1}\omega_i^{-1}  \ ,  \eqno(1\cdot4\rm c)  $$
for $1\leq i\leq N-1 $, where $ \psi_i, \psi_i^\dagger, \omega_i, \omega_i^{-1}$  $(1\leq i\leq N)$ are operators satisfying the following relations:
$$ \omega_i\omega_j=\omega_j\omega_i, \ \ \ \  \omega_i\omega_i^{-1}=\omega_i^{-1}\omega_i=1,   \eqno(1\cdot5\rm a) $$
$$ \omega_i\psi_j\omega_i^{-1}=\psi_j\ \ \ \ (i\neq j),\ \ \  \ \omega_i\psi_i\omega_i^{-1}=q\psi_i, \eqno(1\cdot 5\rm b) $$
$$ \omega_i\psi_j^\dagger\omega_i^{-1}=\psi_j^\dagger \ \ \ \ (i\neq j), \ \ \  \ \omega_i\psi_i^\dagger\omega_i^{-1}=q^{-1}\psi_i^\dagger,  \eqno(1\cdot5\rm c) $$
$$\psi_i\psi_j+\psi_j\psi_i=\psi_i^\dagger\psi_j^\dagger+\psi_j^\dagger\psi_i^\dagger=0,   \eqno(1\cdot5\rm d) $$
$$\psi_i\psi_j^\dagger+\psi_j^\dagger\psi_i=0 \ \ \ (i\neq j) ,  \eqno(1\cdot 5\rm e) $$
$$\psi_i\psi_i^\dagger+q^2\psi_i^\dagger\psi_i=\omega_i^{-2},\ \ \psi_i\psi_i^\dagger+q^{-2}\psi_i^\dagger\psi_i=\omega_i^2.  \eqno(1\cdot 5\rm f) $$ 
Hayashi$^{7)}$denoted the algebra generated by $\psi_i, \psi_i^\dagger, \omega_i,\omega_i^{-1} (1\leq i\leq N)$ by $A_q(N)$ and called it the $q$-Clifford algebra. We should then call $\psi_i,\psi_i^\dagger (1\leq i\leq N) $ the $q$-fermion operators.

\ \ \  As we shall discuss later, however, the $q$-fermion is nothing but the conventional fermion defined by 
$$ \psi_i\psi_j^\dagger+\psi_j^\dagger\psi_i=\delta_{ij},  \eqno(1\cdot 6\rm a)$$ 
$$ \psi_i\psi_j+\psi_j\psi_i=\psi_i^\dagger\psi_j^\dagger+\psi_j^\dagger\psi_i^\dagger=0   \eqno(1\cdot 6\rm b) $$
$$ (1\leq i,j\leq N)   $$
It turns out that all the relations in $(1cdot5)$ are obtained from $(1\cdot6)$ if we define $\omega_i$ and $\omega_i^{-1}$ by 
$$ \omega_i=\psi_i\psi_i^\dagger+q^{-1}\psi_i^\dagger\psi_i,   \eqno(1\cdot7\rm a)$$
$$ \omega_i^{-1}=\psi_i\psi_i^\dagger+q\psi_i^\dagger\psi_i.   \eqno(1\cdot 7\rm b) $$
After cpmpleting  this paper, we became aware of the Note added in the proof of Ref.7), in which we found that the above fact was first pointed out by  Takeuchi but  is yet unpublished.$^{8)}$ We thus see that the algebra $A_q(N)$ can be identified with the conventional fermion algera $A(N)$ defined by $(1\cdot6)$.

\ \ \ We shall discuss  further properties of $A(N)$ and seek the $\ast$-homomorphism $\delta:\ A(N)\longrightarrow A(N)\otimes A(N)$, i.e., a linear mapping from $A(N)$ to $A(N)\otimes A(N)$ satisfying $\delta(ab)=\delta(a)\delta(b),(\delta(a))^*=\delta(a^*)\  a,b,a^*\in A(N)$ and $\delta(1)=1\otimes 1. $
Under certain restrictions, we obtain two $\ast$-homomorphisms, $\delta_1$ and $\delta_2$. Although $\delta_1$ and $\delta_2$ do not satisfy  co-associativity and hence cannot be regarded as the coproduct of $A(N)$, they satisfy a relation similar to  co-associativity. In contrast to $A(N)$, there exists a coproduct for $U_q(A_{N-1})$.  The manner in which the coproduct $\Delta: U_q(A_{N-1})\longrightarrow U_q(A_{N-1})\otimes U_q(A_{N-1})$ should act on $e_i, f_i, k_i, k_i^{-1} (1\leq i\leq N-1)$ is well known. It is given by
$$ \Delta(e_i)=k_i\otimes e_i+e_i\otimes k_i^{-1},  \eqno(1\cdot 8\rm a)$$
$$ \Delta(f_i)=k_i\otimes f_i+f_i\otimes k_i^{-1},  \eqno(1\cdot 8\rm b)$$
$$ \Delta(k_i)=k_i\otimes k_i.  \eqno(1\cdot 8\rm c)  $$
 Since elements of $U_q(A_{N-1})$ and $A(N)$ are related by $(1\cdot4)$, it is expected that the $\Delta$ can be expressed by $\delta_1$ and/or $\delta_2$. We find that this is the case.  

This paper is organized as follows. In $\S 2$, we discuss  how the $q$-fermionic relations $(1\cdot 5)$ are deduced from the conventional fermionic relations $(1\cdot 6)$ with the $\omega_i$ and $\omega_i^{-1}$ defined by $(1\cdot 7)$. In $\S 3$, we consider a mapping from the fermion algebra to its tensor product and obtain two $\ast$-homomorphisms $\delta_1$ and $\delta_2$. We investigate the relation of $\delta_1$ and $ \delta_2$ to the coproduct $\Delta$ of the quantum group $U_q(su(N))$, i.e., $U_q(A_{N-1})$ supplemented with the $\ast$-property mentioned below, in $\S 4$. In $\S 5$, 
we describe a simple application of $\delta_1$ and $\delta_2$. The final section, $\S 6$, is devoted to summary.
\vspace{2.5cm}

\begin{center}
{\large \bf \S 2. ${ q}$-fermion identified with conventional fermion}\par
\end{center}
\vspace{0.5cm}
\ \ \ We consider conventional fermion operators $\psi_i, \psi_i^\dagger(1\leq i\leq N)$ satisfying the anti-commutation relations $(1\cdot 6)$. If we define $ \omega_i $  by $(1\cdot 7\rm a)$, we readily see that its inverse is given by $(1\cdot 7\rm b)$. Moreover, we obtain 
$$ \omega_i^n=\psi_i\psi_i^\dagger+q^{-n}\psi_i^\dagger\psi_i,   \eqno(2\cdot1) $$
for any integer $n$. Here we have made use of the relations $(\psi_i\psi_i^\dagger)^2=\psi_i\psi_i^\dagger,\ (\psi_i^\dagger\psi_i)^2=\psi_i^\dagger\psi_i,\ (\psi_i\psi_i^\dagger)(\psi_i^\dagger\psi_i)=0$ and $(\psi_i^\dagger\psi_i)(\psi_i\psi_i^\dagger)=0$. Putting $n=2$ or $ -2$ in $(2\cdot 1)$, we are led to $(1\cdot5\rm f)$. From the definitions $(1\cdot7)$ of $\omega_i$ and $\omega_i^{-1}$, we get 
$$(\omega_i-1)(q\omega_i-1)=0,   \eqno(2\cdot2\rm a) $$  
$$ \omega_i\psi_i=\psi_i, \ \ \ \psi_i\omega_i=q^{-1}\psi_i,  \eqno(2\cdot2\rm b) $$
$$ \psi_i^\dagger\omega_i=\psi_i^\dagger, \ \ \ \omega_i\psi_i^\dagger=q^{-1}\psi_i^\dagger.   \eqno(2\cdot2\rm c) $$
Equations $(2\cdot2\rm b)$ and $(2\cdot2\rm c)$ should be regarded as the detailed versions of  $(1\cdot 5\rm b)$ and $(1\cdot 5\rm c)$ and the former   reproduces the latter. Thus we see that the conventional fermions reproduce the $q$-fermions.
Although there exists a freedom to define a $q$-fermion as a constant multiple of a conventional fermion, we hereafter fix $ \psi_i,\psi_i^\dagger, \omega_i $ and $ \omega_i^{-1}$ by $(1\cdot6)$ and $(1\cdot7)$ and regard the  $q$-fermions equivalent to the  conventional fermions. Guided by the above discussion, we  consider  the algebra $A(N)$ solely generated by $\psi_i,\psi_i^\dagger(1\leq i\leq N)$ satisfying $(1\cdot6)$. It should be noted that the relation $(2\cdot2\rm a)$ was neither used nor noted by Hayashi. It turns out, however, to be important to discuss the algebra homomorphism of $A(N)$. 

In $A(N)$, we define the $\ast$-operation by 
          $$ (\psi_i)^\ast=\psi_i^\dagger,  \ \ \ (\psi_i^\dagger)^\ast=\psi_i. \eqno(2\cdot3) $$
          $$(c\phi\phi')^\ast=c^*(\phi')^\ast(\phi)^\ast, \ \ \ ( \phi, \ \phi' \in A(N)) \eqno(2\cdot4) $$
where $c$ and $c^*$ are a complex number and its complex conjugate, respectively. Since we are assuming that $q$ is a real number, we have
$$             
(\omega_i)^\ast=\omega_i, \ \ \ (\omega_i^{-1})^\ast=\omega_i^{-1}. \eqno(2\cdot5) $$

The $\ast$-operations for $U_q(A_{N-1})$ are naturally defined with the help of $(2\cdot3),\sim \ (2cdot 5)$ and $(1\cdot 4)$.They are given by 
$$(e_i)^*=f_i, \ \ \ (f_i)^*=e_i, \eqno(2\cdot6\rm a)$$
$$(k_i)^*=k_i, \ \ \ (k_i^{-1})^*=(k_i)^{-1} \eqno(2\cdot6\rm b)$$
The group  $U_q(A_{N-1})$ with the $\ast$-property $(2\cdot6)$ is  denoted by $U_q(su(N))$ and has some  interesting applications in physics. In other words, to consider $U_q(su(N))$, the $\ast$-property $(2\cdot 6)$ is indispensable.

The representation $(1\cdot 4)$ is called the spinor representation of $U_q(su(N))$. Its representation space is the fermion Fock space $V$  spanned by the vectors 
$$ \vert {\bf  m}\rangle\equiv(\psi_1^\dagger)^{m_1}(\psi_2^\dagger)^{m_2}\cdots(\psi_N^\dagger)^{m_N}\vert 0\rangle,  \eqno(2\cdot7)$$
where ${\bf  m}$ is given by
$$ {\bf m}=(m_1,m_2,\cdots,m_n)\in \{0,1\}^N  \eqno(2\cdot8)$$
and $\vert 0\rangle$ is the vacuum satisfying 
$$ \psi_i\vert 0\rangle=0.\ \ \ \ (1\leq i\leq N)  \eqno(2\cdot9)$$
If we denote the vector $(0,0,\cdots,0,1,0,\cdots,0)$ with the $i$th component equal to $1$ and all other components  $0$ by ${\bf  e}$, we have   
$$ \psi_i\vert {\bf m} \rangle=(-1)^{m_1+m_2+\cdots+m_{i-1}}\delta_{m_i,1}\vert {\bf  m}-{\bf e_i}\rangle, \eqno(2\cdot10\rm a)$$
$$ \psi_i^\dagger\vert {\bf  m} \rangle=(-1)^{m_1+m_2+\cdots+m_{i-1}}\delta_{m_{i,0}}\vert {\bf  m}+{\bf e_i}\rangle, \eqno(2\cdot10\rm b)$$
and hence
$$\omega_i\vert {\bf  m}\rangle=q^{-m_i}\vert {\bf  m}\rangle. \eqno(2\cdot10\tt c)$$
The space $V$ is decomposed as 
 $$ V=\displaystyle{\bigoplus_{r=0}^N}V_r,   \eqno(2\cdot11)$$
 where $V_r$ is a subspace of $V$ defined by 
 $$ V_r =\displaystyle{\bigoplus_{\vert{\bf  m}\vert=r}}{\bf  C}\vert {\bf  m}\rangle, \ \ \  \vert {\bf  m} \vert\equiv \displaystyle{\sum_{i=1}^N} m_i. \eqno(2\cdot12)$$
It is known$^{7)}$ that the space $V$ is irreducible under the actions of $A_q(N)$, while each $V_r(0\leq r\leq N)$ is irreducible under  the action of $U_q(A_{N-1})$ with $e_i, f_i, k_i, k_i^{-1}$ given by $(1\cdot4)$.
 
We note that, in the representations mentioned above, the generators $e_i,f_i,k_i,k_i^{-1}$ defined by $(1\cdot4)$ satisfy relations in addition to those in $(1\cdot1)$:
$$ e_i^2=f_i^2=0,$$
$$(k_i-1)(k_i-q)(k_i-q^{-1})=0,$$
$$ e_if_ie_i=e_i,\ \ \ f_ie_if_i=f_i, \ \ \ \ $$
$$k_i=1+\frac{1-q}q(f_ie_i-qe_if_i),$$
$$k_i^{-1}=1+\frac{1-q}q(e_if_i-qf_ie_i),$$  
$$ e_ik_i=q^{-1}e_i,\ \ \ f_ik_i=qf_i,$$
$$ k_ie_i=qe_i,\ \ \ k_if_i=q^{-1}f_i. \eqno(2\cdot13) $$
For example, we have $e_i^2v=0$ for any vector $v$ belonging  to $V$. The  equality $\Delta(e_i^2)=\{\Delta(e_i)\}^2$, however, is in general a non-vanishing operator on $V\otimes V$.
To obtain general representations, we are to consider tensor products of spinor representations.
For operators $a_1,\ a_2,\ \ \cdots,\ a_n,\ \ b_1,
\newline  b_2,\  \cdots,\ b_n\in A(N)$, we define the multplication rule  of tensor products by 
$$ (a_1\otimes a_2\otimes \cdots a_n)(b_1\otimes b_2\otimes\cdots b_n)=(-1)^Ma_1b_1\otimes a_2b_2\otimes\cdots a_nb_n,  \eqno(2\cdot14)$$

$$ M=\displaystyle{\sum_{i=2}^n\sum_{j=1}^{i-1}}d(a_i)d(b_j),  \eqno(2\cdot15)$$
where $d(a)=0$ (respectively $1$) if $a\in A(N)$ consists of an even(respectively, odd) number of fermions operators. The $\ast$-operation on $a_1\otimes a_2\otimes \cdots \otimes a_n$ is defined by 
$$ (a_1\otimes a_2\otimes \cdots\otimes a_n)^*=(-1)^L(a_1)^*\otimes(a_2)^*\otimes\cdots\otimes(a_n)^*, \eqno(2\cdot16)$$
$$ L=\displaystyle{\sum_{i=2}^n\sum_{j=1}^i} d(a_i)d(a_j) .\eqno(2\cdot17)$$
The key tool to discuss the tensor protuct representations of $A(N)$ preserving the $\ast$-property is the $\ast$-homomorphism  $\delta$ of $A(N)$, i.e. a linear mapping from $A(N)$ to $A(N)\otimes A(N)$ satisfying 
$$
\delta(ab)=\delta(a)\delta(b),\ \  \delta(1)=1\otimes 1,\ \ (\delta(a))^*=\delta(a^*),\eqno(2\cdot18) $$

$$(a,b,a^*\in A(N).\ \ \delta(a), \delta(b), \delta(a^*)\in A(N)\otimes A(N))
$$
which we seek in the next section.\par

\vspace{2cm}
\begin{center}

{\large \bf \S 3. Homomorphic mapping of $A(N)$ to $A(N)\otimes A(N)$}\par
\end{center}
\vspace{0.5cm}

We  now  seek a $\ast$-homomorphism $\delta$ of $A(N)$.
We asuume that $\delta(\psi_i)$ takes the form
$$ \delta(\psi_i)=\alpha_i\otimes\psi_i+\psi_i\otimes\beta_i,  \eqno(3\cdot1) $$ 
where $\alpha_i$ and $\beta_i$ consist only of $\psi_i$ and $\psi_i^\dagger$
and satisfy $d(\alpha_i)=d(\beta_i)=0$. Because of the relations $(1\cdot6)$ and $(2\cdot1)$, the most general form of $\alpha_i$ and $\beta_i$ are given by $\alpha_i=a_i\psi_i\psi_i^\dagger+b_i\psi_i^\dagger\psi_i$,$ \beta_i=c_i\psi_i\psi_i^\dagger+d_i\psi_i^\dagger\psi_i,\ a_i,b_i ,c_i,d_i\in {\bf  C}$. Assuming that $a_i,b_i, c_i $ and $ d_i$ are independent of $i$, we have

$$\alpha_i=a\psi_i\psi_i^\dagger+b\psi_i^\dagger\psi_i, $$
$$\beta_i=c\psi_i\psi_i^\dagger+d\psi_i^\dagger\psi_i,  $$
$$a, b, c, d\in{\bf  C}.   \eqno(3\cdot2)$$
We are thus led to consider a $\ast$-homomorphism of the form 
$$ \delta(\psi_i)=(a\psi_i\psi_i^\dagger+b\psi_i^\dagger\psi_i)\otimes\psi_i+\psi_i\otimes(c\psi_i\psi_i^\dagger+d\psi_i^\dagger\psi_i),  \eqno(3\cdot3\rm a) $$
and
$$\delta(\psi_i^\dagger)=(a^*\psi_i\psi_i^\dagger+b^*\psi_i^\dagger\psi_i)\otimes\psi_i^\dagger+\psi_i^\dagger\otimes(c^*\psi_i\psi_i^\dagger+d^*\psi_i^\dagger\psi_i). \eqno(3\cdot3\rm b) $$
The requirements $\delta(\psi_i)\delta(\psi_j)+\delta(\psi_j)\delta(\psi_i)=0$, $\delta(\psi_i^\dagger)\delta(\psi_j^\dagger)+\delta(\psi_j^\dagger)\delta(\psi_i^\dagger)=0$ and  $\delta(\psi_i)\delta(\psi_j^\dagger)+\delta(\psi_j^\dagger)\delta(\psi_i)=0\ \ (i\neq j)$ are obtained from $(1\cdot6)$ and $(2\cdot14)$ irrespectively of the values of $a,b,c,d. $ Another requirement $\delta(\psi_i)\delta(\psi_i^\dagger)+\delta(\psi_i^\dagger)\delta(\psi_i)=1\otimes 1$ is satisfied if $a,b,c$ and $d$ are related by 
$$ \vert a\vert^2+\vert c\vert^2=1$$
$$ \vert a\vert^2=\vert b\vert^2,\ \ \ \vert c\vert^2=\vert d\vert^2,\ \  $$
$$ a^*c-b^* d=0.\eqno(3\cdot4) $$
If we further require that $\delta$ satisfies the co-associativity $(\delta\otimes id)\circ\delta(\phi)=(id\otimes\delta)\circ\delta(\phi),\ \phi\in A(N),\  id=$identity, we are left only with  the trivial homomorphisms $\delta(\phi)=1\otimes\phi$ and $ \delta(\phi)=\phi\otimes1.$ In other words, under the restrictions $(3cdot1)$ and $(3cdot2)$ , there exists no nontrivial coassociative $\ast$-homomorphism for $A(N)$. A somewhat looser requirement 
$$ m((id\otimes\delta)\circ\delta(\phi))=m((\delta\otimes id)\circ\delta(\phi)), \ \ \ (\phi\in A(N))  \ \eqno(3\cdot5) $$
allows  essentially four solutions, where $m$ is the multiplication
$$
m(a_1\otimes a_2\otimes\cdots\otimes a_n)=a_1a_2\cdots a_n.\ \ \ \ \ 
(a_1, a_2,\cdots,a_n\in A(N))  \ \eqno(3\cdot6)$$

Two of the above four solutions satisfy a slightly modified coassociativity:
$$   Y\circ (\delta_k\otimes id)\circ \delta_k\circ Y(a)=(id\otimes \delta_k)\circ\delta_k(a). \ \ (a\in A(N),\ \  k=1,2)   \eqno(3\cdot7) $$

Here the linear mapping $Y$ is defined by 
$$ Y(a_1a_2)=Y(a_2)Y(a_1),   \eqno(3\cdot8\rm a) $$
$$Y(a_1\otimes a_2\otimes \cdots \otimes a_n)=Y(a_n)\otimes\cdots\otimes Y(a_2)\otimes Y(a_1),  \ \ \eqno(3\cdot8\rm b) $$
$$(a_1,a_2,\cdots,a_n\in A(N)) $$
and by 
$$ Y(\psi_i)=\psi_i,\ \ \  Y(\psi_i^\dagger)=\psi_i^\dagger. \ \ (1\leq i\leq N) \eqno(3\cdot9) $$ 
It can be readily seen that $Y$ satisfies 
 $$ Y^2=id. \eqno(3\cdot10) $$
When an algebra homomorphism satisfies the co-associativity, the coalgebra can be defined as an associative algebra. In the case that we  have only $(3\cdot7)$, the coalgebra $B(N)$ dual to $A(N)$ is not associative but turns  out to be a flexible algebra.$^{9)}$ The two homomorphisms satisfying conditions $(3\cdot5)$ as well as $(3\cdot7)$ are given by 
$$ \delta_1(\psi_i)=\frac i{\sqrt 2}(\zeta_i\otimes\psi_i-\psi_i\otimes\zeta_i),$$
$$ \delta_1(\psi_i^\dagger)=-\frac i{\sqrt 2}(\zeta_i\otimes \psi_i^\dagger-\psi_i^\dagger\otimes\zeta_i), \  \eqno(3\cdot11) $$
$$ (1\leq i\leq N)$$
and
$$ \delta_2(\psi_i)=\frac 1{\sqrt 2 }(1\otimes\psi_i +\psi_i\otimes 1),   $$
$$ \delta_2(\psi_i^\dagger)=\frac 1{\sqrt 2}(1\otimes \psi_i^\dagger +\psi_i^\dagger\otimes 1), \eqno(3\cdot12) $$
$$( 1\leq i\leq  N)$$
where $\zeta_i$ is defined by 
$$\zeta_i=\psi_i\psi_i^\dagger-\psi_i^\dagger\psi_i,  \ \ (1\leq i\leq N)  \eqno(3\cdot13)$$
The equality $(3\cdot7)$ is derived in the following way. By direct calculation, we see that Equation $(3\cdot7)$ holds for $a=\psi_i$ and $a=\psi_i^\dagger$. For example,we have 
$$ {\begin{array}{l}
Y\circ (id\otimes \delta_1)\circ\delta_1\circ Y(\psi_i) \\ 
\ \ \ \\
\ \ \ \ =-\frac 1 2 (\psi_i\otimes\zeta_i\otimes\zeta_i-\zeta_i\otimes\psi_i\otimes\zeta_i)+
\frac i{2\sqrt 2}(\zeta_i\otimes 1\otimes\psi_i+1\otimes\zeta_i\otimes\psi_i)  \\
\ \ \\
\ \ \ \ \ \ \ + \frac i{\sqrt 2}(\psi_i^\dagger\otimes\psi_i\otimes\psi_i-\psi_i\otimes\psi_i^\dagger\otimes\psi_i) \end{array}}$$
$$ =(\delta_1\otimes id)\circ\delta_1(\psi_i).
  \eqno(3\cdot14) $$
From the definition of $\delta_1$ and $\delta_2$, we  observe the relation 
$$ d(a)=d(a_{k,l}')+d(a''_{k,l}), \ \ (\rm  mod 2)\ \  \eqno(3\cdot15) $$
where $a'_{k.l} $ and $a''_{k,l}$ are defined by 
$ \delta_k(a)=\sum_l a'_{k,l}\otimes a''_{k,l}.$
 With the aid of (3cdot15), it can be shown that the mappings $U_k, W_k:\  A(N)\to A(N)\otimes A(N) \otimes A(N)$ defined by
$$ U_k=Y\circ(\delta_k\otimes id)\circ\delta_k\circ Y,\ \  W_k=(id\otimes\delta_k)\circ\delta_k   \eqno(3\cdot16) $$
satisfy
$$ 
U_k (ab)=U_k(a)U_k(b),\eqno(3\cdot17\rm a)$$
$$W_k(ab)=W_k(a)W_k(b),
 \eqno(3\cdot17\rm b) $$
for $a,b\in A(N), \ k=1,2.$ This discussion yields $(3\cdot7)$. We have thus obtained two fermionic homomorphisms $ \delta_1 $ and $\delta_2$ by requiring $(3cdot3),(3\cdot5)$ and $(3\cdot7)$.

We here mention that, as for a type of $q$-bosonic algebra, there has been invented a coproduct.$^{10)}$ This coproduct,  however, does not possess,  the $\ast$-property necessary for the discussion of $U_q(su(N))$. In our case of fermion algebra, we respect the $\ast$-property and modify the co-associativity to $(3\cdot7)$.  We  thus  see that the fermion algebra $A(N)$ is a bialgebra with the homomorphism $\delta_1$ or $\delta_2$. The higher dimensional representations of $U_q(su(N))$ can be obtained by making use of $(1\cdot4)$ and $\delta_1$ or $\delta_2$. For example, $\delta_k(\psi_i\psi_{i+1}^\dagger),\ \delta_k(\psi_{i+1}\psi_i^\dagger),\ \delta_k((\psi_i\psi_i^\dagger+q^{-1}\psi_i^\dagger\psi_i)(\psi_{i+1}\psi_{i+1}^\dagger+q\psi_{i+1}^\dagger\psi_{i+1}))$ with $1\leq i\leq N-1,\ k=1$ or $2$ , constitute the representation of $e_i,\ f_i,\ k_i\in U_q(su(N))$, respectively, on $V\otimes V$. The mappings $U_k$ and $W_k$ defined by $(3\cdot16)$ can be used to construct re!
presentions of $U_q(su(N))$ on $V\otimes V\otimes V$.  

We note that $\delta_1$ together with $\delta_2$ generates a $2N$-fermion algebra on $A(N)\otimes A(N).$ This result is obtained by observing
$$ \begin{array}{l}
\delta_1(\psi_i)\delta_2(\psi_j^\dagger)+\delta_2(\psi_j^\dagger)\delta_1(\psi_i)=0,\ \ (i\neq   j) \\
\ \ \ \ \ \\
\delta_2(\psi_i)\delta_1(\psi_j^\dagger)+\delta_1(\psi_j^\dagger)\delta_2(\psi_i)=0,\ \  (i\neq j)  \end{array} \eqno(3\cdot18)$$
and
$$ \delta_1(\psi_i)\delta_2(\psi_j)+\delta_2(\psi_j)\delta_1(\psi_i)=\delta_1(\psi_i^\dagger)\delta_2(\psi_j^\dagger)+\delta_2(\psi_j^\dagger)\delta_1(\psi_i^\dagger)=0 \eqno(3\cdot19) $$
and recalling that both  $\delta_1$ and $\delta_2$ are  homomorphisms of the $N$-fermion algebra. If we define $\Psi_I,\ 1\leq I\leq 2N $ by
$$\begin{array}{lll}
 \Psi_i=\delta_1(\psi_i), &  \Psi_i^\dagger=\delta_1(\psi_i^\dagger),& (1\leq i \leq N) \\
\ & \ & \ \ \\
 \Psi_{i+N} =\delta_2(\psi_i),& \Psi_{i+N}^\dagger=\delta_2(\psi_i), &  (1\leq i\leq N)\end{array} \eqno(3\cdot20)$$
the properties mentioned above yield the $2N$-fermion algebra on $A(N)\otimes A(N)$:
$$\begin{array}{c}
 \Psi_I\Psi_J^\dagger + \Psi_J^\dagger \Psi_I=\delta_{IJ} 1\otimes 1,\\
\ \ \\
 \Psi_I\Psi_J+\Psi_J\Psi_I=\Psi_I^\dagger\Psi_J^\dagger+\Psi_J^\dagger\Psi_I =0.\\
\ \\
(1\leq I,\ \ J\leq 2N)
\end{array} \eqno(3\cdot21)$$ 

\vspace{3.5cm}
\begin{center}
{\large \bf \S4.  $\Delta$ in term of {\bf $\delta_1$ }and {\bf $\delta_2$}}\par
\end{center}
\vspace{0.5cm}
As was discussed in the previous sections, the quantum group $U_q(A_{N-1})$ is related to the algebra $A(N)$ generated by conventional fermion operators. The coproduct $\Delta$ of $U_q(A_{N-1})$ maps $U_q(A_{N-1})$ into $U_q(A_{N-1})\otimes U_q(A_{N-1})$ homomorphically. On the other hand, under centain restrictions we have found two homomorphic mappings $\delta_1,\delta_2:\ A(N)\to A(N)\otimes A(N)$. It would be  very interesting if there existed  some relations between $\Delta$ and $\delta_1$ and/or $\delta_2.$ 

To explore this possibility, it is sufficient to express $\Delta(e_i), \Delta(f_i), \Delta(k_i),\Delta(k_i^{-1}), 1\leq i\leq N-1,$ in terms of $\psi_i,\psi_i^\dagger, 1\leq i\leq N$ and $\delta_k,k=1,2$. We first discuss $\Delta(k_i)$. From $(1\cdot7{\rm a}), (3\cdot11)$ and $(3\cdot12)$, we obtain 
$$ \begin{array}{l}
2\delta _1(\omega_i)=1\otimes\omega_i+\omega_i\otimes 1+{\displaystyle\frac{q-1}q}(\psi_i\otimes \psi_i^\dagger-\psi_i^\dagger\otimes\psi_i),\\
\ \ \ \\
2\delta_2(\omega_i)=1\otimes\omega_i+\omega_i\otimes 1-{\displaystyle\frac{q-1}q}(\psi_i\otimes \psi_i^\dagger-\psi_i^\dagger\otimes\psi). \end{array}\eqno(4\cdot1)$$
Their product is given by $$
\delta _1(\omega_i)\delta_2(\omega_i)= \delta_2(\omega_i)\delta_1(\omega_i)=\omega_i\otimes \omega_i,\ \ \eqno(4\cdot2)   $$
where we have made use of $(2\cdot2\rm a)$. We also have 
$$ \delta_1(\omega_i^{-1})\delta_2(\omega_i^{-1})=\delta_2(\omega_i^{-1})\delta_1(\omega_i^{-1})=\omega_i^{-1}\otimes\omega_i^{-1}.\  \   \eqno(4\cdot3)$$
We here  note the  interesting formula
$$ \delta_1(\omega_i^n)+\delta_2(\omega_i^n)=1\otimes\omega_i^n+\omega_i^n\otimes 1,     \eqno(4\cdot4) $$
where $n$ is an arbitrary integer. Recalling $(1\cdot4\rm c)$ and $(1\cdot8\rm c)$, we readily obtain
$$ \Delta(k_i)=\delta_1(k_i)\delta_2(k_i)=\delta_2(k_i)\delta_1(k_i),   \eqno(4\cdot5\rm a) $$
$$ \Delta(k_i^{-1})=\delta_1(k_i^{-1})\delta_2(k_i^{-1})=\delta_2(k_i^{-1})\delta_1(k_i^{-1}) \eqno(4\cdot5\rm b) $$
from $(4\cdot2)$ and $(4\cdot3)$.
Turning to  discussion of $\Delta(e_i)$ and $\Delta(f_i)$, we define $\tilde \delta$ and a mapping $Z$: $A(N)\otimes A(N)\to A(N)\otimes A(N) $ by
$$  \tilde\delta (a) =\frac1{2\sqrt q}\{(q+1)\delta_2(a)-i(q-1)\delta_1(a)\},\ \  (a\in A(N)) \eqno(4\cdot6) $$
$$  Z(a\otimes b)=Y^{d(a)d(b)}(a\otimes b)=\left\{\begin{array}{lll} a\otimes b & \quad {\rm if} & \quad d(a)d(b)=0, \\
\ \ & \ \ & \ \\
Y(a\otimes b) & \quad {\rm if} & \quad  d(a)d(b)=1.
\end{array}\right.  \eqno(4\cdot7) $$
 Then, with the aid of the relations
$$\omega_i=\frac 1{2}\{(q^{-1}+1)-(q^{-1}-1)\zeta_i\},  \eqno(4\cdot8) $$
$$ \omega_i^{-1}=\frac12\{(q+1)-(q-1)\zeta_i\},  \eqno(4\cdot9) $$
we obtain 
$$\begin{array}{l}
\tilde \delta (\psi_i)=\frac 1{\sqrt 2}(\sqrt q \omega_i\otimes\psi_i+\frac 1{\sqrt q}\psi_i\otimes\omega_i^{-1}),\\
\ \ \\
\tilde \delta(\psi_{i+1}^\dagger)=\frac 1{\sqrt 2}(\sqrt q\psi_{i+1}^\dagger\otimes\omega_{i+1}+\frac 1{\sqrt q}\omega_{i+1}^{-1}\otimes\psi_{i+1}^\dagger) 
\end{array} \eqno(4\cdot10) $$
and hence                     
$$\begin{array}{l}
\tilde \delta(\psi _i)\tilde \delta(\psi_{i+1}^\dagger) \\
\ \ \\
\ \ \ =\frac12\{\omega_i\omega_{i+1}^{-1}\otimes \psi_i\psi_{i+1}^\dagger+\psi_i\psi_{i+1}^\dagger\otimes\omega_i^{-1}\omega_{i+1}\\
\ \ \ \ \ \ \ \\ 
\ \ \ \ \ \ \ \ \ -q\omega_i\psi_{i+1}^\dagger\otimes\psi_i\omega_{i+1}+
q^{-1}\psi_i\omega_{i+1}^{-1}\otimes\omega_i^{-1}\psi_{i+1}^\dagger\}\end{array}$$ $$=\frac12\{k_i\otimes e_i+e_i\otimes k_i^{-1}-q(1-Y)(\omega_i\psi_{i+1}^\dagger\otimes\psi_i\omega_{i+1})\}.   \eqno(4\cdot11) $$
Noting the relations
$$ Y(\omega _i)=q^{-1}\omega _i^{-1},  \ \ \ \  \ Y(\omega_i^{-1})=q\omega_i,   $$
$$ Y(k_i^{-1})=k_i,\ \ \ \ \ Y(k_i)=k_i^{-1},   \eqno(4\cdot12)$$
$$ Y(e_i)=-e_i,$$
we obtain
$$\Delta(e_i)=(id+Z)\circ\{\tilde \delta (\psi_i)\tilde \delta( \psi_{i+1}^\dagger)\}.   \eqno(4\cdot13) $$
Similarly, we have
$$ \Delta(f_i)=(id+Z)\circ\{\tilde\delta(\psi_{i+1})\tilde\delta(\psi_i^\dagger)\}.    \eqno(4\cdot14)$$
It can be seen that the relations $(3\cdot18)$ and $(3\cdot19)$ ensure $\Delta(f_i)=\{\Delta(e_i)\}^\ast$.  We note that equations $(4\cdot13)$ and $(4\cdot14)$ are also written as 
$$ 
 \Delta(e_i)=\{2\tilde \delta(\psi_i)\tilde\delta(\psi_{i+1}^\dagger)\}\bigcap(S\otimes S),\eqno(4\cdot15\rm a) $$
$$
 \Delta(f_i)=\{2\tilde \delta(\psi_{i+1})\tilde\delta(\psi_i^\dagger)\}\bigcap(S\otimes S),  \eqno(4\cdot15\rm b)$$
where $S$ is the linear algebra spanned  by $e_i,f_i,k_i,k_i^{-1}$ of $(1\cdot 4)$, or 
$$\Delta(e_i)=\tilde{\delta}(\psi_i)\tilde{\delta}(\psi_{i+1}^\dagger)+\tau(\{\tilde{\delta}(\psi_i^\dagger)\}^\ast\{\tilde{\delta}(\psi_{i+1})\}^\ast),\eqno(4\cdot16\rm a)$$
$$\Delta(f_i)=\tilde{\delta}(\psi_{i+1})\tilde{\delta}(\psi_i^\dagger)+\tau(\{\tilde{\delta}(\psi_{i+1}^\dagger)\}^\ast\{\tilde{\delta}(\psi_i)\}^\ast), \eqno(4\cdot16\rm b) $$
where $\tau$ is defined by 
$$ \tau(a\otimes b)=b\otimes a.\ \ (a,b\in A(N)) \eqno(4\cdot17)$$
\par

\vspace{2.0cm}
\begin{center}
{\large \bf \S 5. An application of $\delta_1$ and $\delta_2$}
\end{center}
\vspace{0.5cm}

As an application of the homomorphisms $\delta_1$ and $\delta_2$, we consider an eigenvalue problem on the product space $ V\otimes V$, where $V$ in the fermion Fock space discussed in $\S2$. We define the operator $H$ by 
$$ H=\sum_{i,j,=1}^N(a_{ij}\psi_i\otimes \psi_j^\dagger +b_{ij}\psi_i^\dagger \otimes \psi_j).\ \ \  (a_{ij}\in {\bf C}) \eqno(5\cdot1)$$
Reguiring that $H$ satisfies $H^\ast=H$ with the $\ast$-operation defined by $(2\cdot16)$, we have $b_{ij}=-a_{ij}^\ast$. The further requiement $(\rm a)$ $\tau (H)=-H$ yeilds $a_{ij}=a_{ij}^\ast,$ while the requiement $({\rm b)}\ \tau(H)=H$ yeilds $a_{ij}=-a_{ij}^\ast$. Through an appropriate unitary transformation, the matrix $a=(a_{ij})$ is diagonalized, and we are left with 
$$ H_a=\sum_{i,j=1}^N a_{ij}(\psi_i\otimes \psi_j^\dagger-\psi_i^\dagger\otimes\psi_j)=\sum_{l=1}^N \lambda_l(\phi\otimes \phi_l^\dagger-\phi_l^\dagger \otimes \phi_l),\ \ \ (\lambda_l, a_{ij}\in {\bf   R})  \eqno(5\cdot2\rm a) $$
$$ H_b=\sum_{i,j=1}^N a_{ij}(\psi_i\otimes\psi_j^\dagger+\psi_i^\dagger\otimes\psi_j)=\sum_{l=0}^N \sigma _l(\phi_l\otimes \phi_l^\dagger+\phi_l^\dagger\otimes \phi_l)).\ \ (\sigma_l,ia_{ij}\in{\bf R}) \ \ \ \ \eqno(5\cdot2{\rm b})$$
 Here $\phi_l$ is a unitary transformation of $\psi_m$:
 $$ \phi_l=\sum_{M=1}^N u_{lm}\psi_m,\ \ \ \ (u_{lm})=u: {\rm unitary}, \eqno(5\cdot3)$$
It t satisfies the same $N$ fermion algebra as $(1\cdot6)$. The eigenvalue problem of $H_a$ is solved in the follwing way. As we saw at the end of $\S 3$, the set of operators $\{ \Phi_I:\ 1\leq I \leq 2N \ \} $ defined by
 $$ \Phi_i=\delta_1(\phi_i),\ \ \Phi_{i+N}=\delta_2(\phi_i), \ \ (1\leq i\leq N) \eqno(5\cdot4) $$
 constitutes a $2N$ fermion algebra. It is easy to rewrite $H_a$ as 
 $$ H_a=\sum_{i=1}^n \lambda_i(\Phi_i^\dagger\Phi_i-\Phi_{i+N}^\dagger \Phi_{i+N}). \eqno(5\cdot5)$$
  We now see that the vector
  $$\begin{array}{c}
   \vert {\bf   M}\rangle=(\Phi_i^\dagger)^{M_1}(\Phi_2^\dagger)^{M_2}\cdots (\Phi_{2N}^\dagger)^{M_{2n}}(\vert 0\rangle\otimes \vert 0\rangle) \\
   \ \ \ \\
   ({\it  M}_I\in\{ 0,1\},\ \ \ 1\leq I\leq 2N) \end{array}
   \eqno(5\cdot6) $$
  is the eigenvector of $H_a$ belonging to the eigenvalue
  $$ E_{\bf  M}=\sum_{i=1}^n \lambda (M_i-M_{i+N}). \eqno(5\cdot7) $$
  The orthonormality and the completeness of $\{\vert {\bf  M}\rangle,\ {\bf M}\in \{0,\ 1\}^{2N} \} $ in $V\otimes V$ is assured by the $2N$-fermion algebraic property of $\{ \Phi_I: 1\leq I\leq 2N \}.$
  
We note that the eigenvalue problem of $H_a$ can be solved with the aid of the two other homomorphisms mentioned between $(3\cdot5)$ and $(3\cdot6)$. 
\

\vspace{2.5cm}
\begin{center}
{\large \bf \S 6. Summary}
\end{center}
\vspace{0.5cm}

In this paper, we have discussed how Hayashi's $q$-fermion algebra $A_q(N)$ is realized in terms of the conventional fermion algebra $A(N)$. A simple but important observation is that the quantity $\omega_i$ defined by $(1\cdot7\rm a)$ satisfies $(2\cdot1)$. We have found  the algebra homomorphisms $\delta_1,\delta_2: A(N)\to A(N)\otimes A(N)$ satisfying the pseudo-co-associativity, $(3\cdot7)$. We have seen that $\delta_1$ and $\delta_2$ can be used to generate a $2N$-fermion algebra on $A(N)\otimes A(N)$. 
 In terms of $\delta_1$ and $\delta_2$, the $U_q(su(N))$ coproducts $\Delta(k_i)$ and $\Delta(k_i^{-1})$ are given by $(4\cdot 5\rm a)$ and $(4\cdot5\rm b)$, respectively. The combination $\tilde{\delta}$ is convenient for the discussion of $\Delta(e_i)$ 
and $\Delta(f_i)$ as is seen in $(4\cdot13)\sim (4\cdot16)$.
The expression for $\Delta$ is  rather simple but is not simple. Some deeper meaning of them  should be  sought.\par

\vspace{2.0cm}
\begin{center}{\bf Acknowledgements}
 \end{center}
 \par
\vspace{0.8cm}
We are grateful to Professor T. Hayashi and Professor M. Takeuchi for correspondences. We thank Professor S. Hamamoto and Professor T. Kurimoto for their kind interest. One of the authors (M.H.) is grateful to Professor Z. Z. Xin for his hearty hospitality at Liaonig University, China, where a part of this paper was written.  We also thank Mr. Hui-Min Zhang for his help in typesetting.

\par
{\small 
\vspace{2.5cm}
\begin{center}
{\large \bf References}
\end{center}
\par
\vspace{0.8cm}
1) V. G. Drinfel'd, in {\it Proceedings of the International Congress of Mathematicians}  

\ \ \ \ (Berkeley,  CA, USA, 1986), p. 793.

2) M. Jimbo,  Lett. Math. Phys. {\bf 10} (1985), 63.

3) S. L. Woronowicz, Publ. RIMS, Kyoto Univ. {\bf 23} (1987), 117.

4) Z. Chang, Phys. Rep.{\bf 262} (1995), 137.

5) L. C. Biedenharn, J. of Phys. {\bf A22} (1989), L873.

6) A. J. Macfarlane, J. of Phys. {\bf A22} (1989), 4581.

7) T. Hayashi, Commun. Math. Phys. {\bf 127} (1990), 129.

8) M. Takeuchi  and  T. Hayash,  private communications.

9) R. D. Schafer, {\it Introduction to Nonassociative Algebras } (Academic Press,
 
\ \ \ \ New York, 1966).

10) Hong Yan, J. of Phys. {\bf A23}  (1990), L1155.  

}

\end{document}